\documentclass[12pt]{amsart}

 \usepackage{fullpage}

\newcommand{\gve}{\varepsilon}

\newcommand{\ra}{\rightarrow}

\def\lfa#1#2{\mathord{\hbox{\vtop{\kern-5pt\hbox{$#2$}}$\setminus$$#1$}}}

\theoremstyle{plain}
\newtheorem{lma}{Lemma}
\newtheorem{thm}[lma]{Theorem}

\theoremstyle{definition}

\swapnumbers
\def\proof#1{\par\medskip\noindent{\it Proof#1.}}
\def\enprule{\vrule height1mm depth1mm width2mm}
\def\endproof{\discretionary{}{\kern\hsize}{\kern3mm}\llap{\enprule}\par\medskip}
\def\P{\mathbb{P}}
\def\Z{\mathbb{Z}}
\def\N{\mathbb{N}}
\def\vm{\vec{m}}
\def\vn{\vec{n}}
\def\va{\vec{a}}
\def\vb{\vec{b}}

\def\vx{\vec{x}}
\def\vw{\vec{w}}
\def\vf{\vec{f}}
\def\vzr{\vec{0}}

\def\cB{\mathcal{B}}

\def\IP{\hbox{IP}}
\let\ld=\ldots
\def\nasp{\spacefactor=1000}

\begin{document}

\title{The shifted primes and 
the multidimensional Szemer\'edi and polynomial van der Waerden theorems}

\author{Vitaly Bergelson}
\address{Department of Mathematics, The Ohio State University, Columbus, Ohio 43210}
\email{vitaly@math.ohio-state.edu}

\author{Alexander Leibman} 
\address{Department of Mathematics, The Ohio State University, Columbus, Ohio 43210}
\email{leibman@math.ohio-state.edu}

\author{Tamar Ziegler}
\address{Department of Mathematics, Technion, Haifa, Israel 32000}
\email{tamarzr@tx.technion.ac.il} 

\thanks{The first and the third authors are supported by BSF grant No.\nasp\ 2006094.
The first and the second authors are supported by NSF grant DMS-0901106.}


\maketitle
The goal of this short note is to establish new refinements of 
multidimensional Szemer\'edi and polynomial van der Waerden theorems.
Let $\P$ be the set of positive prime integers.
We provide a short derivation of the following statements:

\begin{thm}\label{main1} 
Let $\vm_1, \ld, \vm_k \in \Z^d$
and let $E$ be of positive upper Banach density in $\Z^d$, namely $d^*(E) = \limsup \frac{|E\cap B|}{|B|}>0$, where the limsup is taken over parallelepipeds $B \subset \Z^d$, $B=\prod_{i=1}^d [M_i,N_i]$ with  $\min_{i} |N_i-M_i| \to  \infty$.  
Then the set
\[
R(E)=\bigl\{n\in\N: d^{*}(E \cap (E-n \vm_1) \cap \ld \cap (E-n \vm_k)) >0\bigr\}
\]
has a nonempty intersection with $\P-1$ and with $\P+1$.
\end{thm}

\begin{thm}\label{main2} 
Let $(X,\cB,\mu)$ be a finite measure space
and let $T_1,\ld,T_k$ be pairwise commuting measure preserving transformations of $X$. 
Let $A\in\cB$, $\mu(A)>0$;
then the set
\[
R(A)=\bigl\{n\in\N: \mu(A\cap T_1^{-n}A \cap \ld \cap T_k^{-n}A) >0\bigr\}
\]
has a nonempty intersection with $\P-1$ and with $\P+1$.
\end{thm}

Theorems~\ref{main1} and \ref{main2} are equivalent
via the Furstenberg correspondence principle.
(See, for example, Theorem~6.4.17 in \cite{b} or Theorem~2.1 in \cite{bc}.)
We also remark that for any integer $a\neq\pm1$
one can easily construct counter examples via periodic sets/systems, 
so that Theorems~\ref{main1} and \ref{main2} do not hold true for $\P+a$.

\medskip
$\IP_{r}$ and $\IP_{r}^{*}$ sets (in $\N$) are defined as follows:

\medskip\noindent
{\bf Definition.}
For $r\in\N$,
{\it an $\IP_{r}$ set\/} in $\N$ is a set of the form 
$\{\vn\cdot\vw\}_{\vw\in\{0,1\}^{r}\setminus\{\vzr\}}$, 
where $\vn=(n_1,\ld,n_r)\in\N^{r}$.
A subset of $\N$ is {\it an $\IP_{r}^{*}$ set\/}
if it has a nonempty intersection with every $\IP_{r}$ set in $\N$.
\medskip

For example, an $\IP_{3}$ set is a $7$-element set of the form $\{n,m,k,n+m, n+k,m+k,n+m+k\} \subset \N  $.
\medskip

Our proof of Theorem~\ref{main1} is based on the following two very deep theorems. 
The first was obtained by Furstenberg and Katznelson in \cite{fk} (see Theorem $10.1$ and the remark on page $168$):

\begin{thm}\label{IP}
For any probability measure space $(X,\cB,\mu)$, 
any commuting measure preserving transformations $T_1,\ld,T_k$, 
and any set $A\in\cB$ of positive measure, 
there exists an integer $r$ such that the  set $R(A)$ is an $\IP_{r}^{*}$ set.
\end{thm}

The second was obtained in a series of papers by Green, Tao and Ziegler 
in \cite{gt-primes}, \cite{gt-mobius}, \cite{gtz}.

\begin{thm}\label{GTZ}
Let $\psi_{1},\ld,\psi_{l}$ be affine linear forms in $r$ variables
with integer coefficients,
$\psi_{i}(\vx)=\sum_{j=1}^{r}m_{i,j}x_{j}+c_{i}$, 
no two of which are affinely dependent.
Then there exists an $\vn\in\Z^{r}$ 
such that $\psi_{1}(\vn),\ld,\psi_{l}(\vn)\in\P$ 
iff for any $k\in\N$, $k\geq 2$, there exists an $\vx\in\Z^{r}$
such that $\psi_{1}(\vx),\ld,\psi_{l}(\vx)$ are all nondivizible by $k$.  
\end{thm}

As a corollary, we get that the set $\P-1$ 
(as well as the set $\P+1$)
contains an $\IP_{r}$ set in $\N$ for every $r\in\N$.
Indeed, for any $r$,
since $1$ is not divisible by any $k\geq 2$,
by Theorem~\ref{GTZ} there exists $\vn\in\Z^{r}$ 
such that the integers $\vw\cdot\vn+1$, $\vw\in\{0,1\}^{r}\setminus\{\vzr\}$, 
are all prime.

\proof{ of Theorem~\ref{main2}} 
By Theorem~\ref{IP},
$R(A)$ nontrivially intersects any $\IP_{r}$ set in $\N$ for $r$ large enough,
and by Theorem~\ref{GTZ}, 
the sets $\P-1$ and $\P+1$ contain $\IP_{r}$ sets for all $r$.%
\endproof

We remark that in the case $d=1$ and $T_{i}=T^{i}$, $i=1,\ld,k$,
Theorems~\ref{main1}, \ref{main2} were proved in \cite{fhk} 
conditional on the {\em inverse conjecture for the Gowers norms} 
which was recently obtained in \cite{gtz}. 
However, in their full generality,  
Theorems~\ref{main1}, \ref{main2} cannot be obtained by the methods in that paper. 

We also remark that one cannot obtain 
polynomial extensions of Theorems~\ref{main1} and \ref{main2}
by the methods of the present short note,
since there is so far no polynomial version of Theorem~\ref{IP}. 
(See however, \cite{wz}, where such an extension has been obtained 
for the case $T_{i}$ are all equal to the same transformation.)
On the other hand, 
a ``partition'' version of Theorem~\ref{main1}
(and of the topological version of Theorem~\ref{main2}) 
can be extended to polynomials: 

\begin{thm}\label{main3} 
For any $d\in\N$
and any finite partition $\Z^{d}=\bigcup_{s=1}^{c}C_{s}$,
at least one of the sets $C_{s}$ has the property that
for any finite set of polynomials $\vf_{i}\colon\Z\ra\Z^{d}$, $i=1,\ld,k$,
satisfying $\vf_{i}(0)=0$ for all $i$,
there exist $p\in\P$ and $\va\in\Z^{d}$ such that
$$
\va,\va+\vf_{1}(p-1),\ld,\va+\vf_{k}(p-1)\in C_{s},
$$
and there exist $q\in\P$ and $\vb\in\Z^{d}$ such that
$$
\vb,\vb+\vf_{1}(q+1),\ld,\vb+\vf_{k}(q+1)\in C_{s}.
$$
\end{thm}

A parallel topological dynamical result is the following refinement of Theorem $C$ in \cite{bl1}: 

\begin{thm}\label{main4}
Let $(X,\rho)$ be a compact metric space
and let $T(\vm)$, $\vm\in\Z^{d}$, 
be an action of $\Z^{d}$ on $X$ by continuous transformations.
Then for any finite set of polynomials $\vf_{i}\colon\Z\ra\Z^{d}$, $i=1,\ld,k$, 
with $\vf_{i}(0)=0$ for all $i$,
and any $\gve>0$
there exist a point $x\in X$ and a prime integer $p\in\P$ such that
$\rho\bigl(x,T(\vf_{i}(p-1))x\bigr)<\gve$ for all $i=1,\ld,k$,
and there exist a point $y\in X$ and a prime integer $q\in\P$ such that
$\rho\bigl(y,T(\vf_{i}(q+1))y\bigr)<\gve$ for all $i=1,\ld,k$.
\end{thm}

(See \cite{fw} and \cite{r} for a discussion of equivalence
of Ramsey-theoretical and topological-dynamical results.)

The proof of Theorem~\ref{main3} is the same as of Theorem~\ref{main1},
based on the following version of the polynomial van der Waerden theorem:

\begin{thm}\label{pvdW} 
For any partition $\Z^{d}=\bigcup_{s=1}^{c}C_{s}$
at least one of the sets $C_{s}$ has the property
that for any finite set of polynomials $\vf_{i}\colon\Z\ra\Z^{d}$, $i=1,\ld,k$,
with $\vf_{i}(0)=0$ for all $i$,
$$
\bigl\{n\in\N:
\hbox{$\va,\va+\vf_{1}(n),\ld,\va+\vf_{k}(n)\in C_{s}$ for some $\va\in\Z^d$}\bigr\}
$$
is an $\IP_{r}^{*}$ set for $r$ large enough. 
\end{thm}

\medskip\noindent{\bf Remark.}
The polynomial van der Waerden theorem, proved in \cite{bl1},
was formulated in a slightly weaker form:
it was only claimed there (see \cite{bl1}, Corollary~1.12) that the set 
$$
\bigl\{n\in\N:
\hbox{$\va,\va+\vf_{1}(n),\ld,\va+\vf_{k}(n)\in C_{s}$ 
for some $\va\in\Z^d$ and some $s$}\bigr\}
$$ 
is an $\IP^{*}$ set.
({\it An $\IP^{*}$ set\/} in $\N$ is a set 
that has a nonempty intersection with every $\IP$ set,
where {\it an $\IP$ set\/} is an $\IP_{\infty}$ set,
that is, a set of the form 
$\cup_{k=1}^{\infty} \{ \vn \cdot \vw\}_{\vw\in\{0,1\}^{k}\setminus\{\vzr\}}$
for some $\vn=(n_{1},n_{2},\ld)\in\N^{\N}$.)
However, Theorem~\ref{pvdW} can be easily derived from the results of \cite{bl1}.
Namely, the proof of Corollary~1.9 in \cite{bl1} actually shows 
that the set $P$ in the formulation of this corollary is an $\IP_{r}^{*}$ set 
for $r$ large enough,
and a standard application of this corollary 
to a minimal closed shift-invariant subset of the space of $c$-partitions of $\Z^{d}$
gives the desired result.
Another way to get it 
is by utilizing the polynomial Hales-Jewett theorem (\cite{bl2}).



\begin{thebibliography}{10}

\bibitem{b} V.Bergelson,
{\emph Ergodic theory and Diophantine problems}.
Topics in symbolic dynamics and applications (Temuco, 1997),  167-205, 
London Math. Soc. Lecture Note Ser. {\bf 279}, 
Cambridge Univ. Press, Cambridge, 2000.

\bibitem{bc} V.Bergelson, R. McCutcheon,
{\emph R. Recurrence for semigroup actions and a non-commutative Schur theorem.}  
Topological dynamics and applications (Minneapolis, MN, 1995),  205--222, 
Contemp. Math. {\bf 215}, Amer. Math. Soc., Providence, RI, 1998.

\bibitem{bl1} V. Bergelson, A. Leibman, 
\emph{Polynomial van der Waerden and Szemer\'{e}di theorems}, 
J. of AMS {\bf9} (1996), 725-753.

\bibitem{bl2} V. Bergelson, A. Leibman, 
\emph{Set-polynomials and polynomial extension of the Hales-Jewett theorem}, 
Ann. of Math. (2) {\bf150} (1999), no. 1, 33-75. 

\bibitem{fhk} N. Frantzikinakis, B. Host and B. Kra, 
\emph{Multiple recurrence and convergence for sequences related to the prime numbers}, 
J. Reine Angew. Math. {\bf611} (2007), 131-144.

\bibitem{fk} H. Furstenberg, Y. Katznelson, 
\emph{ An ergodic Szemer?di theorem for IP-systems and combinatorial theory}, 
J. Analyse Math. {\bf45} (1985), 117-168. 


\bibitem{fw} H. Furstenberg, B. Weiss, 
\emph{Topological dynamics and combinatorial number theory},
J. d'Analyse Math. {\bf34} (1978), 61-85.

\bibitem{gt-primes} B. Green, T. Tao, 
\emph{Linear equations in primes}, to appear in Ann. Math.

\bibitem{gt-mobius} B. Green, T. Tao, 
\emph{The M\"obius function is strongly orthogonal to nilsequences}, 
to appear in Ann. of Math.

\bibitem{gtz}
B. Green, T. Tao, T. Ziegler, 
\emph{An inverse theorem for the Gowers $U^{s+1}[N]$ norm}, 
preprint.

\bibitem{r} R. McCutcheon, 
\emph {Elemental methods in ergodic Ramsey theory}. 
Lecture Notes in Mathematics, {\bf 1722}. Springer-Verlag, Berlin, 1999.

\bibitem{wz} T. Wooley, T. Ziegler, 
\emph{Multiple recurrence and convergence along the primes}, 
submitted.

\end{thebibliography}
\end{document}